\documentclass[11pt]{article}

\usepackage{latexsym,amssymb}

\makeatletter \@addtoreset{equation}{section}

\newtheorem{thm}{Theorem}[section]

\newtheorem{lem}[thm]{Lemma}
\newtheorem{cor}[thm]{Corollary}

\def\enu#1{\newline\makebox[5mm][l]{\rm(#1)}}

\def\bp{\noindent{\it Proof.}\ }
\def\bpp#1{\smallskip\noindent{\it Proof of #1.}\ }
\def\ep{\nopagebreak\newline\mbox{\ }\hfill\rule{2mm}{2mm}}
\def\epp{\nopagebreak\mbox{\ }\hfill\rule{2mm}{2mm}}

\def\Ad{{\rm Ad}}

\def\End{{\rm End}}
\renewcommand{\Im}{{\rm Im}}
\def\ind{{\rm Ind}_F\,}
\def\indq{q{\rm-Ind}_F\,}
\def\Hom{{\rm Hom}}

\def\Ker{{\rm Ker}}
\def\Mat{{\rm Mat}}

\def\Tr{{\rm Tr}}

\def\7#1{{\mathbb #1}}

\def\A{{\cal A}}
\def\B{{\cal B}}
\def\C{{\cal C}}

\def\H{{\cal H}}
\def\L{{\cal L}}

\def\LL{{L^2(A,\varphi)}}

\def\eps{{\varepsilon}}

\def\rad{{\blacktriangleleft}}
\def\r{{\triangleleft}}
\def\l{{\triangleright}}

\def\2{{\frac{1}{2}}}

\def\>{{\rangle}}
\def\<{{\langle}}

\def\D{{\Delta}}

\def\Ir0{{\rm Irr}(A_0,\Delta_0)}
\def\sl2{{U_q({\mathfrak su}_2)}}

\begin{document}

\title{\bf Hopf Algebra Equivariant Cyclic Cohomology,
K-theory and Index Formulas}

\author{Sergey Neshveyev$^1$
\& Lars Tuset}

\footnotetext[1]{Supported by the Norwegian Research Council.}

\date{}

\maketitle

\begin{abstract}
For an algebra $\B$ with an action of a Hopf algebra $\H$ we
establish the pairing between equivariant cyclic cohomology
and equivariant $K$-theory for $\B$. We then extend this formalism
to compact quantum group actions and show that equivariant cyclic
cohomology is a target space for the equivariant Chern character
of equivariant summable Fredholm modules. We prove an analogue of
Julg's theorem relating equivariant $K$-theory to ordinary
$K$-theory of the C$^*$-algebra crossed product, and characterize
equivariant vector bundles on quantum homogeneous spaces.
\end{abstract}

\section*{Introduction}

This paper is part of a project which aims to extend the theory of
noncommutative geometry of Connes~\cite{Co} to $q$-homogeneous
spaces in such a way that the various $q$-deformed geometric
notions play as natural a role as possible. This is beneficial
from the point of view of computations. In this paper we focus on
the pairing between equivariant cyclic cohomology and equivariant
$K$-theory, acknowledging the Chern character as the backbone of
noncommutative geometry.

Equivariant cyclic cohomology for group actions was introduced
in~\cite{KKL}. It is straightforward to generalize these results
to the context of Hopf algebras with involutive antipode. In the
general case the pertinent formulas are less obvious. Several
approaches have been proposed in~\cite{AK1,AK2}. In the previous
version of this paper we introduced equivariant cyclic cohomology
in the case when the Hopf algebra admits a modular element, by
which we mean a group-like
element implementing the square of the antipode.
However, this cyclic cohomology turns out to be
isomorphic to the one introduced in the latest version
of~\cite{AK1}, so we prefer working with the latter definition.

In Section~\ref{s2} we establish basic properties of equivariant
cyclic cohomology such as Morita invariance and the construction
of the trace map. In Section~\ref{s3} we set up the pairing
between equivariant cyclic cohomology and the equivariant
$K$-groups. Restricting to the case when the Hopf algebra in
question has a modular element, we discuss the relation of this
pairing to twisted cyclic cohomology~\cite{KMT}. In
Section~\ref{s4} we extend these results to compact quantum groups
in that we replace the Hopf algebra by the algebra of finitely
supported functions on the dual discrete quantum group.
Incorporating equivariant summable even Fredholm modules allows us
to state an index formula. In the remaining part of the paper we
study equivariant $K$-theory for compact quantum group actions. We
establish its relation to ordinary $K$-theory of the crossed
product algebra, thus generalizing a result of Julg~\cite{J}.
Finally we obtain an analogue of the fact that equivariant vector
bundles on homogeneous spaces are induced from representations of
the stabilizer group.

\section{Cyclic cohomology} \label{s2}

\subsection{Cocyclic objects}

Recall \cite{L} that a cocyclic object in an abelian category $\C$
consists of a sequence of objects $C^n$ in $\C$, $n\ge0$, together
with morphisms
$$d^n_i\colon C^{n-1}\to C^n ,\ \ \
s^n_i\colon C^{n+1}\to C^n ,\ \ \ \ t_n\colon C^n\to C^n ,\ \ \ \
(0\le i\le n)$$ that satisfy
\begin{eqnarray*}
d^n_jd^{n-1}_i&=&d^n_id^{n-1}_{j-1}\ \ \hbox{for}\ \ i<j;\\
s^n_js^{n+1}_i&=&s^n_is^{n+1}_{j+1}\ \ \hbox{for}\ \ i\le j;\\
s^n_jd^{n+1}_i&=&d^n_is^{n-1}_{j-1}\ \ \hbox{for}\ \ i<j;\\
s^n_jd^{n+1}_i&=&\iota\ \ \hbox{for}\ \ i=j\ \ \hbox{or}\ \ i=j+1;\\
s^n_jd^{n+1}_i&=&d^n_{i-1}s^{n-1}_j\ \ \hbox{for}\ \ i>j+1;\\
t_nd^n_i&=&d^n_{i-1}t_{n-1}\ \ \hbox{for}\ \ i\ge 1,\ \
\hbox{and}\ \ t_nd^n_0=d^n_n;\\
t_ns^n_i&=&s^n_{i-1}t_{n+1}\ \ \hbox{for}\ \ i\ge 1,\ \
\hbox{and}\ \ t_ns^n_0=s^n_nt^2_{n+1};\\
t^{n+1}_n&=&\iota.
\end{eqnarray*}
Every cocyclic object gives a bicomplex
\[\begin{array}{ccccccccc} \vdots & & \ \vdots &  & \vdots & &
\ \vdots &  & \\  b\uparrow \quad & & -b'\uparrow \quad &  &
b\uparrow \quad &  & -b'\uparrow \quad & &
\\   C^2 & \stackrel{{1-\lambda}}\to &   C^2 &
\stackrel{{  N}}\to &  C^2 & \stackrel{1-\lambda}\to &  C^2 &
\stackrel{  N}\to &\cdots\\  b\uparrow \quad & & - b'\uparrow
\quad &  &  b\uparrow \quad & & -b'\uparrow \quad & &
\\  C^1 & \stackrel{1-\lambda}\to &   C^1 &
\stackrel{ N}\to &  C^1 & \stackrel{1-\lambda}\to &  C^1 &
\stackrel{ N}\to &\cdots\\  b\uparrow \quad & & -b'\uparrow \quad
& & b\uparrow \quad &  & -b'\uparrow \quad & &
\\  C^0 & \stackrel{1-\lambda}\to &   C^0 &
\stackrel{ N}\to &  C^0 & \stackrel{1-\lambda}\to &  C^0 &
\stackrel{ N}\to &\cdots\\
\end{array}\]
with morphisms $b_n=\sum^n_{i=0}(-1)^id^n_i$,
$b'_n=\sum^{n-1}_{i=0}(-1)^id^n_i$, $\lambda_n=(-1)^nt_n$ and
$N^n=\sum^n_{i=0}\lambda^i_n$. The cyclic cohomology $HC^\bullet$
of a cocyclic object is by definition the cohomology of the total
complex associated to its bicomplex, whereas the Hochschild
cohomology $HH^\bullet$ is defined to be the cohomology of the
complex $(C^\bullet,b)$. A cocyclic object yields a long exact
sequence, referred to as the $IBS$-sequence of Connes, written
$$
\dots\longrightarrow HC^n\stackrel{I}{\longrightarrow}HH^n
\stackrel{B}{\longrightarrow}HC^{n-1}\stackrel{S}{\longrightarrow}
HC^{n+1}\stackrel{I}{\longrightarrow}\dots,
$$
where $S$ is the periodicity operator and
$B_n=N^ns^n_nt_{n+1}(\iota-\lambda_{n+1})$.

We denote the cohomology of the complex $(\Ker(\iota-\lambda),b)$
by $H^\bullet_\lambda$. The equality
$HC^\bullet=H^\bullet_\lambda$ holds whenever the objects $C^n$
are vector spaces over a ground field that contains the rational
numbers. In this case the periodicity
operator gives the periodic cyclic cohomology $HP^0$ and $HP^1$ as
direct limits of $HC^{2n}$ and $HC^{2n+1}$, respectively.

\subsection{Equivariant cyclic cohomology}

By $(\H,\D)$ we mean a Hopf algebra (over $\7C$) with
comultiplication $\D$, invertible antipode $S$ and counit
$\varepsilon$. We adapt the Sweedler notation
$$
\D^n(\omega)=\omega_{(0)}\otimes\dots\otimes\omega_{(n)},
$$
(with summation understood) where $\D^n=(\iota\otimes\D)\D^{n-1}$.
Let $\B$ denote a unital right $\H$-module algebra, that is, a
unital algebra with a right action $\r$ of $\H$ on $\B$ such that
$$
(b_1b_2)\r\omega=(b_1\r\omega_{(0)})(b_2\r\omega_{(1)})\ \
\hbox{and}\ \ 1\r\omega=\eps(\omega)1.
$$
Following \cite{AK1} we introduce a cocyclic object in the
category of vector spaces using the above data. The associated
cyclic cohomology is called the {\it equivariant cyclic
cohomology} of $\B$ and denoted by $HC^\bullet_{\H}(\B)$. Likewise
we write $HH^\bullet_{\H}(\B)$ and $HP^\bullet_{\H}(\B)$ for the
{\it equivariant Hochschild cohomology} and the {\it equivariant
periodic cyclic cohomology}, respectively. Let $C^n=C^n_{\H}(\B)$
be the space of linear functionals $f$ on
$\H\otimes\B^{\otimes(n+1)}$ that are $\H$-invariant in the sense
that
$$
\omega\l f=\varepsilon(\omega)f\ \ \hbox{for all}\ \ \omega\in\H.
$$
Here the left action $\l$ of $\H$ on the space of linear
functionals on $\H\otimes\B^{\otimes(n+1)}$
is given by
$$
(\omega\l f)(\eta\otimes b_0\otimes\dots\otimes b_n)
=f(S^{-1}(\omega_{(0)})\eta\omega_{(1)}\otimes
b_0\r\omega_{(2)}\otimes \dots\otimes b_n\r \omega_{(n+2)}).
$$
To see that this actually is a left action, recall that the tensor
product $X_1\otimes X_2$ of two right $\H$-modules $X_1$ and $X_2$
is again a right $\H$-module with action
$$
(x_1\otimes x_2)\r\omega=x_1\r\omega_{(0)}\otimes x_2\r\omega_{(1)}.
$$
Thus, considering $\H$ itself as a right $\H$-module with respect
to the action
$$
\eta\r\omega=S^{-1}(\omega_{(0)})\eta\omega_{(1)},
$$
we regard $\H\otimes\B^{\otimes(n+1)}$ as a right $\H$-module,
which in turn shows that the space of linear functionals on
$\H\otimes\B^{\otimes(n+1)}$ is a left $\H$-module. Now define
\begin{eqnarray*}
(t_nf)(\omega\otimes b_0\otimes\dots\otimes b_n)
&=&f(\omega_{(0)}\otimes b_n\r\omega_{(1)}\otimes
    b_0\otimes\dots\otimes b_{n-1}),\\
(s^n_if)(\omega\otimes b_0\otimes\dots\otimes b_n)
&=&f(\omega\otimes b_0\otimes\dots\otimes b_i\otimes1\otimes
    b_{i+1}\otimes\dots\otimes b_n),\\
(d^n_if)(\omega\otimes b_0\otimes\dots\otimes b_n)
&=&f(\omega\otimes b_0\otimes\dots\otimes
    b_{i-1}\otimes b_ib_{i+1}\otimes b_{i+2}
    \otimes\dots\otimes b_n)\ \ \hbox{for}\ \ i<n
\end{eqnarray*}
and put $d^n_n=t_nd^n_0$. These formulas give linear maps
$$
t_n\colon C^n_{\H}(\B)\to C^n_{\H}(\B),\ \ d^n_i\colon
C^{n-1}_{\H}(\B)\to C^n_{\H}(\B), \ \ s^n_i\colon
C^{n+1}_{\H}(\B)\to C^n_{\H}(\B)
$$
which define a cocyclic object in the category of vector spaces.
Indeed, if we consider a new Hopf algebra $(\H_1,\D_1)$ with
$\H_1=\H$ and $\D_1=\D^{op}$, and define a left action of $\H_1$
on $\B$ by $\omega\l b=b\r S(\omega)$, then the correspondence
$f\mapsto\tilde f$, where
$$
\tilde f(b_0\otimes\dots\otimes b_n)(\omega)=f(S^2(\omega)\otimes
b_0\otimes\dots\otimes b_n),
$$
defines an isomorphism of the cocyclic object above with the
cocyclic object associated with the left $\H_1$-algebra~$\B$ as
defined by~\cite{AK1}. The axioms of the cocyclic objects can also
be checked directly. We include a proof of the most critical
identity $t_n^{n+1}=\iota$:

Note that
$$
(t^{n+1}_nf)(\omega\otimes b_0\otimes\dots\otimes b_n)
=f(\omega_{(0)}\otimes b_0\r\omega_{(1)}\otimes\dots\otimes
b_n\r\omega_{(n+1)}).
$$
Hence it suffices to show that if $X$ is a right $\H$-module and
$f$ is a linear functional on $\H\otimes X$ such that $\eta\l
f=\varepsilon(\eta)f$, then
$$
f(\omega\otimes x)=f(\omega_{(0)}\otimes x\r\omega_{(1)}).
$$
Due to the identity $\omega=\omega_{(0)}\r\omega_{(1)}$ we indeed
obtain
$$
f(\omega_{(0)}\otimes x\r\omega_{(1)})
=f(\omega_{(0)}\r\omega_{(1)}\otimes x \r\omega_{(2)})
=\eps(\omega_{(1)})f(\omega_{(0)}\otimes x)=f(\omega\otimes x).
$$

\subsection{Morita invariance} \label{s2.3}

Let $X$ be a right $\H$-module and let $\End(X)$ denote the
algebra of all complex linear maps on $X$. Write $\pi_X$ for the
antihomomorphism $\H\to\End(X)$ given by the right $\H$-action.
Then $\End(X)$ is clearly a right $\H$-module algebra with respect
to the adjoint action
$$
T\rad\omega=\pi_X(\omega_{(0)})T\pi_XS^{-1}(\omega_{(1)}).
$$
If in addition $X$ is a unital algebra, we identify $X$ with the
subalgebra of $\End(X)$ consisting of those endomorphisms which
are given by multiplication from the left. In this case the
adjoint action on $\End(X)$ restricts to the original action of
$\H$ on $X$ if and only if $X$ is a right $\H$-module algebra.

\begin{lem} \label{2.2}
For any right $\H$-module $X$,
the subalgebra $\End(X)\otimes \B$ of $\End(X\otimes \B)$ is an
$\H$-submodule with respect to the adjoint action of $\H$ on
$\End(X\otimes \B)$, where the tensor product $X\otimes\B$ is
regarded as a right $\H$-module. Moreover, the adjoint action on
$\End(X)\otimes \B$ is given by
$$
(T\otimes b)\rad\omega=\pi_X(\omega_{(0)})T\pi_XS^{-1}
(\omega_{(2)})\otimes b\r\omega_{(1)},
$$
for which $\End(X)\otimes\B$ is a right $\H$-module algebra.
\end{lem}

\bp First observe that
$\pi_{X\otimes\B}(\omega)
=\pi_X(\omega_{(0)})\otimes\pi_\B(\omega_{(1)})$.
To any $b\in\B$ define $L_b$ to be the endomorphism on $\B$ given
by left multiplication with $b$. Then for $T\in\End(X)$, we derive
\begin{eqnarray*}
(T\otimes L_b)\rad\omega
&=&\pi_{X\otimes\B}(\omega_{(0)})(T\otimes L_b)\pi_{X\otimes\B}
S^{-1}(\omega_{(1)})\\
&=&\pi_X(\omega_{(0)})T\pi_XS^{-1}(\omega_{(3)})\otimes
\pi_\B(\omega_{(1)})L_b\pi_\B S^{-1}(\omega_{(2)})\\
&=&\pi_X(\omega_{(0)})T\pi_XS^{-1}(\omega_{(2)})\otimes
L_{b\r\omega_{(1)}}.
\end{eqnarray*}
It is now immediate that $\End(X)\otimes\B$ is an $\H$-submodule
subalgebra of $\End(X\otimes \B)$. \ep

Note that $\B$ is in general not an $\H$-submodule subalgebra of
$\End(X)\otimes\B$, while $\End(X)$ is. Remark also that since
$\End(X)$ and $\B$ are right $\H$-modules, one can consider
the $\H$-module tensor product action
$$
(T\otimes b)\r\omega=T\rad\omega_{(0)}\otimes b\r\omega_{(1)}.
$$
When the Hopf algebra $\H$ is cocommutative this action coincides
with the action above. In general, however, these two actions are
different, and $\End(X)\otimes\B$ equipped with the tensor product
action need not even be an $\H$-module algebra.

Suppose now that $X$ is finite dimensional. For each $n\ge0$ and
$f\in C^n_{\H}(\B)$, define
$$
(\Psi^n f)(\omega\otimes(T_0\otimes b_0)\otimes\dots
\otimes(T_n\otimes b_n)) =f(\omega_{(0)}\otimes
b_0\otimes\dots\otimes b_n) \Tr(\pi_XS^{-1}(\omega_{(1)})T_0\dots
T_n),
$$
where $\Tr$ is the canonical non-normalized trace on $\End(X)$.
\begin{lem}
\label{2.3}
The formula above gives maps
$$
\Psi^n\colon C^n_{\H}(\B)\to C^n_{\H}(\End(X)\otimes\B),
$$
which constitute a morphism of cocyclic objects, and thus induces
homomorphisms
$$[\Psi^n]\colon HC^n_{\H}(\B)\to
HC^n_{\H}(\End(X)\otimes\B).$$
\end{lem}

\bp Let $f\in C^n_{\H}(\B)$. To check $\H$-invariance of
$\Psi^nf$, compute
\newline $(\Psi^n f)((\eta\otimes(T_0\otimes
b_0)\otimes\dots\otimes(T_n\otimes b_n))\r\omega)$
\begin{eqnarray*}
&=&(\Psi^n f)(S^{-1}(\omega_{(0)})\eta\omega_{(1)}\otimes
(\pi_X(\omega_{(2)})T_0\pi_XS^{-1}(\omega_{(4)})
\otimes b_0\r\omega_{(3)})\otimes\dots\\
& &\hspace{1cm}\dots\otimes(\pi_X(\omega_{(3n+2)})T_n
\pi_XS^{-1}(\omega_{(3n+4)})
\otimes b_n\r\omega_{(3n+3)}))\\
&=&f(S^{-1}(\omega_{(1)})\eta_{(0)}\omega_{(2)}\otimes
b_0\r\omega_{(5)}\otimes\dots\otimes b_n\r\omega_{(3n+5)})
\Tr(\pi_XS^{-2}(\omega_{(0)})\pi_XS^{-1}(\eta_{(1)})\times\\
& &\hspace{1cm}\times\pi_XS^{-1}(\omega_{(3)})
\pi_X(\omega_{(4)})T_0\pi_XS^{-1}(\omega_{(6)})
\dots\pi_X(\omega_{(3n+4)})T_n\pi_XS^{-1}(\omega_{(3n+6)})).
\end{eqnarray*}
Using the identity
$\pi_XS^{-1}(\omega_{(0)})\pi_X(\omega_{(1)})=\eps(\omega)$
repeatedly, this expression simplifies to \newline $\displaystyle
f(S^{-1}(\omega_{(1)})\eta_{(0)}\omega_{(2)}\otimes
b_0\r\omega_{(3)}\otimes\dots\otimes b_n\r\omega_{(n+3)})
\Tr(\pi_XS^{-2}(\omega_{(0)})\pi_XS^{-1}(\eta_{(1)})T_0\dots
T_n\pi_XS^{-1}(\omega_{(n+4)}))$
\begin{eqnarray*}
&=&\eps(\omega_{(1)})f(\eta_{(0)}\otimes b_0\otimes\dots\otimes
b_n) \Tr(\pi_XS^{-2}(\omega_{(0)})\pi_XS^{-1}(\eta_{(1)})T_0\dots
T_n\pi_XS^{-1}(\omega_{(2)}))\\
&=&f(\eta_{(0)}\otimes b_0\otimes\dots\otimes b_n)
\Tr(\pi_XS^{-2}(\omega_{(0)}S(\omega_{(1)}))\pi_XS^{-1}
(\eta_{(1)})T_0\dots T_n)\\
&=&\eps(\omega) (\Psi^n f)(\eta\otimes(T_0\otimes b_0)\otimes\dots
\otimes(T_n\otimes b_n)).
\end{eqnarray*}
Thus $\Psi^n$ is a well-defined map $C^n_{\H}(\B)\to
C^n_{\H}(\End(X)\otimes\B)$. It is readily verified that
$\Psi^\bullet$ commutes with the maps $t_n$, $d^n_i$ and $s^n_i$.
\ep

Suppose $p\in\End(X)$ is an idempotent such that
$p\rad\omega=\eps(\omega)p$ and $pX$ is a submodule with trivial
$\H$-action. In other words,
$\pi_X(\omega)p=p\pi_X(\omega)=\eps(\omega)p$. Then the map
$\Phi_p\colon\B\to\End(X)\otimes\B$ given by $\Phi_p(b)=p\otimes
b$ is a homomorphism of algebras which respects the right actions
of $\H$. Hence it induces homomorphisms $[\Phi^n_p]\colon
HC^n_{\H}(\End(X)\otimes\B)\to HC^n_{\H}(\B)$.

\begin{thm} \label{Morita}
Suppose $X$ is a finite dimensional vector space endowed with the
trivial right action $x\r\omega=\eps(\omega)x$ of $\H$. Then
$[\Psi^n]\colon HC^n_{\H}(\B)\to HC^n_{\H}(\End(X)\otimes\B)$ from
Lemma \ref{2.3} is an isomorphism with inverse $[\Phi^n_p]$, where
$p$ is any one-dimensional idempotent in $\End(X)$.
\end{thm}

\bpp{Theorem \ref{Morita}} The demonstration is similar to the
proof in the non-equivariant case: Recall~\cite[Corollary
2.2.3]{L} that the $IBS$-sequence of Connes implies that if a
morphism of cocyclic objects induces an isomorphism for Hochschild
cohomology then it also does so for cyclic cohomology. Since
$\Phi_p\circ\Psi=\iota$ already on the level of complexes, it is
therefore enough to show that $\Psi\circ\Phi_p$ induces the
identity map on Hochschild cohomology. Thus it suffices to show
that $\Psi\circ\Phi_p$ is homotopic to the identity regarded as a
morphism of $(C^\bullet_{\H}(\End(X)\otimes\B),b)$. Let
$x_1,\dots,x_m$ be a basis in $X$ with corresponding matrix units
$m_{ij}\in\End(X)$. We may assume that $p=m_{11}$. Then
\cite[Theorem 1.2.4]{L} the required homotopy $h^n\colon
C^{n+1}_{\H}(\End(X)\otimes\B)\to C^n_{\H}(\End(X)\otimes\B)$ is
given by $h^n=\sum^n_{j=0}(-1)^jh^n_j$, where
$$
(h^n_jf)(\omega\otimes(T_0\otimes
b_0)\otimes\dots\otimes(T_n\otimes b_n))=\sum_{k_0,\dots,k_{j+1}}
(T_0)_{k_0k_1}(T_1)_{k_1k_2}\dots(T_j)_{k_jk_{j+1}}\times
$$
$$
\times f(\omega\otimes(m_{k_01}\otimes b_0)\otimes(m_{11}\otimes
b_1)\otimes\dots\otimes(m_{11}\otimes
b_j)\otimes(m_{1k_{j+1}}\otimes1)\otimes(T_{j+1}\otimes
b_{j+1})\otimes\dots\otimes(T_n\otimes b_n))
$$
and $T_{ij}$ stand for matrix coefficients of $T$.
\ep

The following corollary is a standard consequence of Morita
invariance. It is obtained by considering two embeddings $b\mapsto
m_{11}\otimes b$, $b\mapsto m_{22}\otimes b$ of $\B$ into
$\Mat_2(\7C)\otimes\B=\Mat_2(\B)$ and the automorphism
$\Ad\left(\matrix{1 & 0\cr 0 & b}\right)$ of $\Mat_2(\B)$.

\begin{cor} \label{2.5}
Let $b$ be an invertible $\H$-invariant element of $\B$. Then the
automorphism $\Ad\,b$ induces the identity map on
$HC^\bullet_{\H}(\B)$. \epp
\end{cor}

\bigskip

\section{K-theory} \label{s3}

We make the same assumptions on $(\H,\D)$ and $\B$ as in the
previous section.

\subsection{Equivariant K-theory} \label{s3.1}

By an $\H$-equivariant
$\B$-module we mean a vector space with $\H$-module and
$\B$-module structures compatible in the sense that
$$
(xb)\r\omega=(x\r\omega_{(0)})(b\r\omega_{(1)}).
$$
Equivalently, equivariant modules can be considered
as $(\B\rtimes\H)$-modules, where the crossed product
$\B\rtimes\H$ is the algebra with underlying space $\B\otimes\H$
and product
$$
(b\otimes\omega)(c\otimes\eta)=b(c\r
S^{-1}(\omega_{(1)}))\otimes\omega_{(0)}\eta.
$$
In the sequel we will always write $b\omega$ instead of
$b\otimes\omega$ whenever we mean an element of $\B\rtimes\H$. The
identities
$$
\omega b=(b\r S^{-1}(\omega_{(1)}))\omega_{(0)}\ \ \hbox{and}\ \
b\omega=\omega_{(0)}(b\r\omega_{(1)})
$$
in $\B\rtimes\H$ are easily verified.

Let $\C_\H(\B)$ be the full subcategory of the category of
$(\B\rtimes\H)$-modules consisting of those modules which are
isomorphic to $X\otimes\B$ for finite dimensional (as vector
spaces) right $\H$-modules~$X$. As before, we consider
$X\otimes\B$ as a right $\H$-module with tensor product action,
while the $\B$-module structure is defined by the right action of
$\B$ on itself. Furthermore, let $\tilde\C_\H(\B)$ denote the
associated pseudo-abelian category, that is, the category of
modules isomorphic to $X_p=p(X\otimes\B)$ for idempotents
$p\in\End_{\B\rtimes\H}(X\otimes\B)$. Note that for
$T\in\End_\B(X\otimes\B)=\End(X)\otimes\B$ we have
$T\in\End_\H(X\otimes\B)$ if and only if $T$ is $\H$-invariant.
Thus $\End_{\B\rtimes\H}(X\otimes\B)$ is the algebra of
$\H$-invariant elements in $\End(X)\otimes\B$.

Let $K_0^\H(\B)$ be the Grothendieck group of the category
$\tilde\C_\H(\B)$. More concretely, consider all possible
$\H$-invariant idempotents $p\in\End(X)\otimes\B$. Say that two
such idempotents $p\in\End(X)\otimes\B$ and
$p'\in\End(X')\otimes\B$ are equivalent if there exist
$\H$-invariant elements $\gamma\in\Hom_\7C(X,X')\otimes\B$ and
$\gamma'\in\Hom_\7C(X',X)\otimes\B$ such that $\gamma\gamma'=p'$
and $\gamma'\gamma=p$. Note that $\Hom_\7C(X,X')\otimes\B$ can be
considered as a subspace of $\End(X\oplus X')\otimes\B$ and as
such is endowed with a canonical $\H$-module structure. The set of
equivalence classes of $\H$-invariant idempotents is an abelian
semigroup in a standard fashion. Then $K^\H_0(\B)$ coincides with
the associated Grothendieck group of this semigroup.

Recall also that if for two $\H$-invariant idempotents
$p\in\End(X)\otimes\B$ and $p'\in\End(X')\otimes\B$ there exist
$\H$-invariant elements $\gamma\in\Hom_\7C(X,X')\otimes\B$ and
$\gamma'\in\Hom_\7C(X',X)\otimes\B$ such that $\gamma\gamma'=p'$
and $\gamma'\gamma=p$, then, if we consider $p$ and $p'$ as
elements of $\End(X\oplus X')\otimes\B$, there exists an
$\H$-invariant invertible element $\gamma_0\in\End(X\oplus
X')\otimes\B$ such that $\gamma_0p\gamma^{-1}_0=p'$. Namely, put
$$
\gamma_0=\pmatrix{1-p & p\gamma'p'\cr p'\gamma p & 1-p'}.
$$
Hence we can equivalently define $K^\H_0(\B)$ as the Grothendieck
group of the semigroup of equivalence classes of $\H$-invariant
idempotents with the equivalence relation generated by similarity
(two idempotents $p$ and $p'$ are called similar if there exists
an invertible $\H$-invariant element
$\gamma\in\Hom_\7C(X,X')\otimes\B$ such that $\gamma
p\gamma^{-1}=p'$) and the condition $p\sim p\oplus0$, where we
think of $p\oplus0$ as living in $\End(X\oplus X')\otimes\B$.

\smallskip

To define $K_1^\H(\B)$, consider invertible elements
$u\in\End_{\B\rtimes\H}(Y)$, where $Y$ is an object in
the category $\C_\H(\B)$. Then $K_1^\H(\B)$ is by
definition the abelian group generated by isomorphism
classes of such $u$ satisfying the following relations:
\enu{i} $[u]+[v]=[uv]$ for $u,v\in\End_{\B\rtimes\H}(Y)$;
\enu{ii} $[u_1]+[u_2]=[u]$ whenever there exists a {\it split}
exact sequence
$$
0\longrightarrow Y_1\stackrel{i}{\longrightarrow}
Y\stackrel{j}{\longrightarrow}Y_2\longrightarrow0
$$
with $u_1\in\End_{\B\rtimes\H}(Y_1)$,
$u_2\in\End_{\B\rtimes\H}(Y_2)$ and $u\in\End_{\B\rtimes\H}(Y)$
such that $ui=iu_1$ and $u_2j=ju$.

The relation (ii) can be rewritten as follows. Given invertible
$\H$-invariant elements $u_i\in\End(X_i)\otimes\B$, $i=1,2$, and
an $\H$-invariant element $T\in\Hom_\7C(X_2,X_1)\otimes\B$, we
have
$$
[u_1]+[u_2]=\left[\pmatrix{u_1 & T\cr 0 & u_2}\right].
$$

Note also that if we use $\tilde\C_\H(\B)$ instead of $\C_\H(\B)$
in the definition of $K_1^\H(\B)$, we get exactly the same group.

\smallskip

If $\H$ is semisimple, which in particular implies that it is
finite dimensional, we can equivalently describe the
category $\tilde\C_\H(\B)$ as the category of f.g.
(finitely generated) projective $(\B\rtimes\H)$-modules.

\begin{thm} \label{Julg1}
Suppose $\H$ is semisimple. Then there exist bijective
correspondences between \enu{i} the equivalence classes of
$\H$-invariant idempotents in $\End(X)\otimes\B$, where $X$ ranges
over all finite dimensional $\H$-modules; \enu{ii} the isomorphism
classes of $\H$-equivariant f.g. projective right $\B$-modules;
\enu{iii} the isomorphism classes of f.g. projective right
$(\B\rtimes\H)$-modules.

In particular, $K_0^\H(\B)\cong K_0(\B\rtimes\H)$
and $K_1^\H(\B)\cong K_1(\B\rtimes\H)$.
\end{thm}

\bp For any $\H$-invariant idempotent $p\in\End(X)\otimes\B$ we
have already seen that $X_p=p(X\otimes\B)$ is an $\H$-equivariant
f.g. projective right $\B$-module. It is furthermore clear that
two idempotents $p\in\End(X)\otimes\B$ and
$p'\in\End(X')\otimes\B$ are equivalent if and only if $X_p$ and
$X'_{p'}$ are isomorphic as $\H$-equivariant $\B$-modules, or as
$(\B\rtimes\H)$-modules. So to prove the theorem we just have to
show that \enu{a} $X_p$ is a projective $(\B\rtimes\H)$-module;
\enu{b} any $\H$-equivariant f.g. projective right $\B$-module is
isomorphic to $X_p$ for some $X$ and $p$; \enu{c} any f.g.
projective right $(\B\rtimes\H)$-module is isomorphic to $X_p$ for
some $X$ and $p$.

Since $X_p$ is a direct summand of $X\otimes\B$, in proving (a) we
may assume $p=1$. By semisimplicity $X$ is isomorphic to a finite
direct sum of modules $q_1\H,\dots,q_n\H$ for some idempotents
$q_1,\dots,q_n\in\H$. Thus we may assume $X=q\H$ for some
idempotent $q\in\H$. The map $q\H\otimes\B\to q(\B\rtimes\H)$,
$\omega\otimes b\mapsto \omega b$, is an isomorphism with inverse
obtained by restricting the map $\B\rtimes\H\to\H\otimes\B$,
$b\omega\mapsto\omega_{(0)}\otimes b\r\omega_{(1)}$, to
$q(\B\rtimes\H)$. Thus $X_p$ is indeed a f.g. projective
$(\B\rtimes\H)$-module.

To establish (b) consider an $\H$-equivariant f.g. projective
right $\B$-module $Y$ with generators $y_1,\dots,y_n$. Set
$X=y_1\r\H+\dots+y_n\r\H$. Then $X$ is an $\H$-module which is
finite dimensional as a vector space. The map $T\colon
X\otimes\B\to Y$, $x\otimes b\mapsto xb$, is a surjective
homomorphism of $(\B\rtimes\H)$-modules. Since $Y$ is a projective
$\B$-module, there exists a $\B$-module map $T'\colon Y\to
X\otimes\B$ such that $TT'=\iota$. By semisimplicity there exists
a right integral $\eta$ in $\H$ (that is,
an element such that $\eta\omega=\eps(\omega)\eta$ for any
$\omega\in\H$)  such that $\eps(\eta)=1$.
Replacing $T'$ by $T'\rad\eta=\pi_{X\otimes
\B}(\eta_{(0)})T'\pi_YS^{-1}(\eta_{(1)})$, it is clear that we can
assume that $T'$ is also an $\H$-module map. Then $Y\cong X_p$
with $p=T'T$. The proof of (c) is similar, but shorter since $T'$
from the beginning can be chosen to be a $(\B\rtimes\H)$-module
map.

The established bijective correspondences immediately imply
$K_0^\H(\B)\cong K_0(\B\rtimes\H)$.
Since every exact sequence of projective modules is split,
the result for $K_1$ also follows.
\ep

\subsection{Pairing with cyclic cohomology}

Let $R(\H)$ be the space of $\H$-invariant linear functionals on
$\H$. As before, we consider the action
$\eta\r\omega=S^{-1}(\omega_{(0)})\eta\omega_{(1)}$ of $\H$ on
itself.

\begin{thm}
There exists pairings
$$
\<\cdot,\cdot\>\colon HC^{2n}_{\H}(\B)\times K^\H_0(\B) \to R(\H)\ \
\hbox{and}\ \
\<\cdot,\cdot\>\colon HC^{2n+1}_{\H}(\B)\times K^\H_1(\B) \to R(\H),
$$
such that for $f\in C^{2n}_{\H}(\B)$ and $p\in\End(X)\otimes\B$ we
have
$$
\<[f],[p]\>(\omega)=(\Psi^{2n}f)(\omega\otimes
p\otimes\dots\otimes p),
$$
and such that for $f\in C^{2n+1}_{\H}(\B)$ and
$u\in\End(X)\otimes\B$ we have
$$
\<[f],[u]\>(\omega)=(\Psi^{2n+1}f)(\omega\otimes(u^{-1}-1)\otimes
(u-1)\otimes\dots\otimes(u^{-1}-1)\otimes(u-1)),
$$
where $\Psi^\bullet$ is the map defined before Lemma~\ref{2.3}.
\end{thm}

\bp The proof can be reduced to the non-equivariant case as
follows. First observe that the quantities
$\<f,p\>(\omega)=(\Psi^{2n}f)(\omega\otimes p\otimes\dots\otimes
p)$ and
$$
\<f,u\>(\omega)=(\Psi^{2n+1}f)(\omega\otimes(u^{-1}-1)\otimes
(u-1)\otimes\dots\otimes(u^{-1}-1)\otimes(u-1))
$$
remain unchanged if we replace $p$ by $p\oplus0$ and $u$ by
$u\oplus1$. Secondly note that the relation
$$
\<f,u_1\>(\omega)+\<f,u_2\>(\omega)
=\<f,\pmatrix{u_1 & T\cr 0 & u_2}\>(\omega)
$$
follows from $\<f,uv\>(\omega)=\<f,u\>(\omega)+\<f,v\>(\omega)$,
because $\pmatrix{u_1 & T\cr 0 & u_2}=\pmatrix{u_1 & 0\cr 0 & u_2}
\pmatrix{1 & u_1^{-1}T\cr 0 & 1}$ and $\<f,u\>(\omega)=0$, whenever
$(u-1)^2=0$, due to the condition $t_{2n+1}f=-f$.

By substituting $\B$ with $\End(X)\otimes\B$ it is thus enough to
show that $\<f,p\>(\omega)$ and $\<f,u\>(\omega)$ depend only on
the cohomology class of $f$, the similarity class of
$p\in\B^\H$ and the isomorphism class of $u\in\B^\H$, and that
$\<f,uv\>(\omega)=\<f,u\>(\omega)+\<f,v\>(\omega)$ for invertible
elements $u$ and $v$ in $\B^\H$. Now for a fixed $\omega\in\H$
consider the morphism $\omega_*\colon C^\bullet_\H(\B)\to
C^\bullet(\B^\H)$ of cocyclic objects given by
$$
(\omega_*f)(b_0\otimes\dots\otimes b_n)
=f(\omega\otimes b_0\otimes\dots\otimes b_n)
$$
for $f\in C^n_\H(\B)$. Then $\<f,p\>(\omega)=\<\omega_*f,p\>$ and
$\<f,u\>(\omega)=\<\omega_*f,u\>$, so the result follows from the
analogous result in the non-equivariant case \cite{Co}. \ep

Note that the proof in \cite{Co} implies even a stronger result
that $HC^{2n+1}_\H(\B)$ pairs with $K_1^\H(\B)/\sim_h$, where
$\sim_h$ is the equivalence relation given by polynomial homotopy.

\smallskip

Suppose now that we have a group-like element $\rho$ of $\H$.
Consider the twisted cyclic cohomology
$HC^\bullet_{\theta_\rho}(\B)$ of $\B$ introduced in \cite{KMT},
where $\theta_\rho$ is the twist automorphism of $\B$ given by
$\theta_\rho (b)=b\r\rho$. More generally, for any finite
dimensional right $\H$-module $X$ consider the twist automorphism
$\theta_{X,\rho}=\Ad\,\pi_X(\rho)\otimes\theta_\rho$ on $\End
(X)\otimes\B$. As in Subsection \ref{s2.3} we can construct a map
$$
\Psi^n_\rho\colon HC^n_{\theta_\rho}(\B)\to
HC^ n_{\theta_{X,\rho}}(\End (X)\otimes\B)
$$
given by
$$
(\Psi^n_\rho f)((T_0\otimes b_0)\otimes\dots
\otimes(T_n\otimes b_n)) =f(b_0\otimes\dots\otimes b_n)
\Tr(\pi_X(\rho^{-1})T_0\dots T_n).
$$
Using this map we  get a pairing
$$
\<\cdot,\cdot\>_\rho\colon
HC^{2n}_{\theta_\rho}(\B)\times K^\H_0(\B) \to \7C
$$
given by $\<[f],[p]\>_\rho =(\Psi_\rho^{2n}f)(p^{\otimes(2n+1)})$.
On the other hand, there exists a map $\rho_*\colon
HC^{n}_\H(\B)\to HC^{n}_{\theta_\rho}(\B)$ determined by
$$(\rho_*\phi)(b_0\otimes\dots\otimes b_n)
=\phi(\rho\otimes b_0\otimes\dots\otimes b_n),$$
and such that the following diagram
\[\begin{array}{ccc}
HC^{2n}_\H(\B)\times K^\H_0(\B) & \stackrel{\<\cdot,\cdot\>}
\longrightarrow & R(\H)\\
\rho_* \downarrow\quad &  & \rho_*\downarrow\quad\\
HC^{2n}_{\theta_\rho}(\B)\times K^\H_0(\B) &
\stackrel{\<\cdot,\cdot\>_\rho}
\longrightarrow & \7C
\end{array}\]
commutes, where $\rho_*\colon R(\H)\to\7C$ is the evaluation map
$\rho_*(f)=f(\rho)$. This diagram is particularly interesting for
Hopf algebras which admit a modular element $\rho$.
We shall return to it in Subsection~\ref{s4.2}.

\bigskip

\section{Compact quantum group actions} \label{s4}

\subsection{Extending the formalism}

So far we have only considered Hopf algebra actions. The theory
does however apply to more general contexts, notably to the case
when $\H$ is the algebra of finitely supported functions on a
discrete quantum group. As explained in Introduction this is in
fact the main motivation for the present work, and  the purpose of
this section is to extend our results to this setting.

Throughout this section  $(A,\D)$ denotes a compact quantum group
in the sense of Woronowicz, and $\alpha\colon B\to A\otimes B$
denotes a left coaction of $(A,\D)$ on a unital C$^*$-algebra $B$.
We shall follow notation and conventions in \cite{NT} with the
only exception that we replace $\rho$ by $\rho^{-1}$. In
particular, we assume that $(A,\D)$ is a reduced compact quantum
group, so the Haar state $\varphi$ on $A$ is faithful. The
coaction $\alpha$ is assumed to be non-degenerate in the sense
that $(A\otimes1)\alpha(B)$ is dense in $A\otimes B$.

The dual discrete quantum group $(\hat A,\hat\D)$ has a canonical
dense $*$-subalgebra $\hat\A$, which can be considered as the
algebra of finitely supported functions on the discrete quantum
group. The elements of $\hat\A$ are bounded linear functionals on
$A$. Thus the coaction $\alpha$ defines a right action of $\hat\A$
on $B$ given by
$$
b\r\omega=(\omega\otimes\iota)\alpha(b).
$$
Set $\B=B\r\hat\A$. The assumption on non-degeneracy ensures that
$\B$ is a dense $*$-subalgebra of $B$ such that
$\alpha(\B)\subset\A\otimes\B$ and $\B\r\hat\A=\B$, where
$\A\subset A$ is the Hopf $*$-algebra of matrix coefficients of
irreducible unitary corepresentations of $(A,\D)$.

The $*$-algebra $\hat\A$ is isomorphic to the algebraic direct sum
$\oplus_{s\in I}B(H_s)$ of full matrix algebras $B(H_s)$, where
$I$ denotes the set of equivalence classes of irreducible unitary
corepresentations of $(A,\D)$. The comultiplication $\hat\D$
restricts to a $*$-homomorphism
$$
\hat\D\colon\hat\A\to M(\hat\A\otimes\hat\A)
=\prod_{s,t\in I}B(H_s)\otimes B(H_t).
$$
If $\hat\D(\omega)=(\omega_{st})_{s,t\in I}$, then the matrix
$(\omega_{st})_{s,t}$ has only finitely many non-zero components
in each column and row. This together with the property
$\B\r\hat\A=\B$ shows that expressions like
$b_1\r\omega_{(0)}\otimes b_2\r\omega_{(1)}$ make sense as
elements of $\B\otimes\B$. It is then straightforward to check
that all the results and notions from the previous sections remain
true for $\B$ whenever $(\H,\D)$ is replaced by $(\hat\A,\hat\D)$.
Even Theorem~\ref{Julg1}, which is not directly applicable since
$\hat\A$ is not semisimple, has an analogue in this context, as
will be shown in Subsection~\ref{s4.3} below.

Remark that if $(A,\D)$ is a $q$-deformation of a compact
semisimple Lie group $G$, then $\A'=M(\hat\A)$ has a dense Hopf
subalgebra, namely the $q$-deformed universal enveloping algebra
of the Lie algebra of the group $G$. One can work with this
algebra instead of $\hat\A$. Then our previous results do not
require any modifications. Since this algebra consists of
unbounded elements, in many cases it is however more convenient to
work with $(\hat\A,\hat\D)$.

\subsection{The Chern character and index formulas} \label{s4.2}

Suppose $(\pi,H,F,\gamma,U)$ is an $\alpha$-equivariant even
Fredholm module over $\B$. By this we mean that $\pi$ is a bounded
$*$-representation of $\B$ on a graded Hilbert space $H=H_-\oplus
H_+$, that $F$ is a symmetry on $H$ ($F=F^*$, $F^2=1$), that
$\gamma=\pmatrix{-1 & 0\cr 0 & 1}$ is the grading operator for
which $F$ is odd and $\pi(b)$, $b\in\B$, is even, and that the
commutator $[F,\pi(b)]$ is a compact operator on $H$ for each
$b\in\B$. Furthermore, the element $U\in M(A\otimes K(H))$ is a
unitary corepresentations of $(A,\D)$, so
$(\D\otimes\iota)(U)=U_{13}U_{23}$, which commutes with $1\otimes
F$ and with $1\otimes\gamma$, and
$$
U(\iota\otimes\pi)\alpha(b)=(1\otimes\pi(b))U.
$$
The corepresentation $U$ defines a bounded $*$-representation of
$\hat\A$ on $H$ by $\pi_U(\omega)=(\omega\otimes\iota)(U)$. Recall
that a Fredholm module $(\pi,H,F)$ over $\B$ is called
$(2n+1)$-summable if $[F,\pi(b)]\in\L^{2n+1}(H)$ for any $b\in\B$,
where $\L^{2n+1}(H)$ is the Schatten $(2n+1)$-class. Then we can
define
$$
\phi_F(\omega\otimes b_0\otimes\dots\otimes b_{2n})
=\frac{(-1)^n}{2}
\Tr(\gamma\pi_U(\omega)F[F,\pi(b_0)]\dots [F,\pi(b_{2n})]).
$$
The identity $U(\iota\otimes\pi)\alpha(b)=(1\otimes\pi(b))U$
implies
$\pi_U(\omega_{(0)})\pi(b\r\omega_{(1)})=\pi(b)\pi_U(\omega)$. It
is now easy to see that $\phi_F$ is a cyclic cocycle and thus
defines an element of $HC^{2n}_{\hat\A}(\B)$.

Suppose $p\in\B$ is a $\hat\A$-invariant projection, and define
projections $p_+=\2(1+\gamma)\pi(p)$ and $p_-=\2(1-\gamma)\pi(p)$.
Then $p_-Fp_+\colon p_+H_+\to p_-H_-$ is a Fredholm operator.
Since $F$, $\gamma$ and $\pi(p)$ all commute with the algebra
$\pi_U(\hat\A)$, the finite dimensional spaces $\Ker(p_-Fp_+)$ and
$\Im(p_-Fp_+)^\perp$ are $\hat\A$-invariant. Denote by
$\tilde\pi_+$ and $\tilde\pi_-$ the finite dimensional
representations of $\hat\A$ on $\Ker(p_-Fp_+)$ and
$\Im(p_-Fp_+)^\perp$, respectively, given by restriction. Define
$\phi_\pm(\omega)=\Tr(\tilde\pi_\pm(\omega))$ and $(\ind
p)(\omega)=\phi_+(\omega)-\phi_-(\omega)$. Note that if $\omega$
is a sufficiently large projection in $\hat\A$, then
$$
(\ind p)(\omega)=\dim\tilde\pi_+ - \dim\tilde\pi_-.
$$

\begin{thm}
We have $\displaystyle
(\ind p)(\omega)=\<[\phi_F],[p]\>(\omega)$.
\end{thm}

\bp The proof is analogous to the one in \cite{Co}, and is based
on the fact that if $\tau$ is trace on an algebra, then
$\tau((p-p')^{2n+1})=\tau(p-p')$ for any two idempotents $p$ and
$p'$ and any $n\in\7N$. \ep

Consider now a finite dimensional unitary corepresentation $V\in
A\otimes B(H_V)$ of $(A,\D)$. Set $\tilde U=U_{13}V_{12}\in
M(A\otimes K(H_V\otimes H))$ and consider the coaction
$\tilde\alpha_V$ of $(A,\D)$ on $B(H_V)\otimes B$ given by
$$
\tilde\alpha_V(T\otimes b)=(V^*\otimes1)(1\otimes T\otimes
1)\alpha(b)_{13}(V\otimes1).
$$
Then $(\iota\otimes\pi,H_V\otimes H,1\otimes F,1\otimes
\gamma,\tilde U)$ is an $\tilde\alpha_V$-equivariant even Fredholm
module over $B(H_V)\otimes\B$. For any $\tilde\alpha_V$-invariant
projection $p\in B(H_V)\otimes\B$, we define $\ind p$ as above.
Note that the right action of $\hat\A$ on $B(H_V)\otimes\B$
corresponding to $\tilde\alpha_V$ in the sense that
$x\r\omega=(\omega\otimes\iota)\tilde\alpha_V(x)$ is precisely the
action introduced in Lemma~\ref{2.2}, where  $H_V$ is the right
$\hat\A$-module with action $\xi\r\omega=\pi_V\hat S(\omega)\xi$.
Indeed,
\begin{eqnarray*}
(\omega\otimes\iota)\tilde\alpha_V(T\otimes b)
&=&(\omega\otimes\iota\otimes\iota)((V^*\otimes1)
(1\otimes T\otimes 1)\alpha(b)_{13}(V\otimes1))\\
&=&((\omega_{(0)}\otimes\iota)(V^*)\otimes1)(T\otimes
b\r\omega_{(1)})
((\omega_{(2)}\otimes\iota)(V)\otimes1)\\
&=&\pi_V\hat S(\omega_{(0)})T\pi_V(\omega_{(2)})\otimes
b\r\omega_{(1)}.
\end{eqnarray*}
Thus $p$ defines an element of $K_0^{\hat\A}(\B)$ and we get a
homomorphism $\ind\colon K_0^{\hat\A}(\B)\to R(\hat\A)$. The
cocycle corresponding to the Fredholm module
$(\iota\otimes\pi,H_V\otimes H,1\otimes F,1\otimes \gamma,\tilde
U)$ is exactly $\Psi^{2n}\phi_F$. Hence $\ind=\<[\phi_F],\cdot\>$
on $K_0^{\hat\A}(\B)$.

Recall now that there exists a canonical element $\rho\in
M(\hat\A)$ such that $\hat\D(\rho)=\rho\otimes\rho$ and $\hat
S^2(\omega)=\rho^{-1}\omega\rho$. Then in the above notation
$(\ind p)(\rho)=\phi_+(\rho)-\phi_-(\rho)$ is nothing but the
difference between the quantum dimensions of $\tilde\pi_+$ and
$\tilde\pi_-$. We denote this quantity by $\indq p$ and think of
it as a quantum Fredholm index of $p_-Fp_+$.

We say that an $\alpha$-equivariant Fredholm module
$(\pi,H,F,\gamma,U)$ over $\B$ is $(2n+1,\rho)$-summable if
$$
\pi_U(\rho)^{1\over 2(2n+1)}|[F,\pi(b)]|\pi_U(\rho)^{1\over
2(2n+1)}\in\L^{2n+1}(H)
$$
for any $b\in\B$. Under this condition
$\indq=\<\tilde\phi_F,\cdot\>_\rho$, where $\tilde\phi_F\in
HC^{2n}_{\theta_\rho}(\B)$ is the twisted cocycle given by
$$
\tilde\phi_F(b_0\otimes\dots\otimes b_{2n}) =\frac{(-1)^n}{2}
\Tr(\gamma
F[F,\pi(b_0)]\dots [F,\pi(b_{2n})]\pi_U(\rho)).
$$
As will be shown elsewhere, the quantum index is much easier to
compute than the usual one and the equivariant index $\ind$.

\subsection{Julg's theorem} \label{s4.3}

We have shown how an $\alpha$-equivariant Fredholm module gives
rise to an index map on equivariant $K$-theory. In the rest of the
paper we address the question of how to compute the $K$-theory.
First we shall obtain an analogue of Theorem~\ref{Julg1} and thus
extend the result of~\cite{J} to compact quantum groups.

Equivariant $K$-theory for coactions of arbitrary locally compact
quantum groups was defined in~\cite{BS}. While the general case
requires $KK$-theory technique, the case of compact quantum groups
can be handled using the simple-minded definition in
Subsection~\ref{s3.1}. Thus we define $K^{\hat\A}_0(B)$ as the
Grothendieck group associated with the semigroup of equivalence
classes of all $\hat\A$-invariant projections in $B(H_V)\otimes B$
for all finite dimensional unitary corepresentations of $V$ of
$(A,\D)$ (in view of~\cite{BS} this group should rather be denoted
by $K^A_0(B)$). Since $(B(H_V)\otimes
B)\rad\hat\A=B(H_V)\otimes\B$, by definition
$K^{\hat\A}_0(B)=K^{\hat\A}_0(\B)$. We shall now describe
$K^{\hat\A}_0(B)$ in terms of $\alpha$-equivariant Hilbert
$B$-modules.

By an $\alpha$-equivariant Hilbert $B$-module we mean a right
Hilbert $B$-module $X$ together with a non-degenerate continuous
linear map
$\delta\colon X\to A\otimes X$ (here we consider $A\otimes X$ as a
right Hilbert $(A\otimes B)$-module, and non-degeneracy means that
$\delta(X)(A\otimes1)$ is dense in $A\otimes X$) such that
$(\D\otimes\iota)\delta=(\iota\otimes\delta)\delta$ and
\enu{i}
$\delta(\xi b)=\delta(\xi)\alpha(b)$ for $\xi\in X$ and $b\in B$;
\enu{ii} $\<\delta(\xi_1),\delta(\xi_2)\>_{A\otimes B}
=\alpha(\<\xi_1,\xi_2\>_B)$ for $\xi_1,\xi_2\in X$.

If $V\in A\otimes B(H_V)$ is a finite dimensional unitary
corepresentation of $(A,\D)$ and $p\in B(H_V)\otimes B$ is a
$\hat\A$-invariant projection, then $H_{V,p}=p(H_V\otimes B)$ is
an $\alpha$-equivariant Hilbert $B$-module with
$$
\delta(\xi\otimes
b)=(V^*\otimes1)(1\otimes\xi\otimes1)\alpha(b)_{13}\ \ \hbox{and}\
\ \<\xi_1\otimes b_1,\xi_2\otimes b_2\>=(\xi_2,\xi_1)b_1^*b_2.
$$

\begin{lem}
Any $\alpha$-equivariant f.g. Hilbert $B$-module $X$ is isomorphic
to $H_{V,p}$ for some $V$ and $p$.
\end{lem}

\bp The coaction $\delta$ defines a right action of $\hat\A$ on
$X$, $\xi\r\omega=(\omega\otimes\iota)\delta(\xi)$. By
non-degeneracy $X\r\hat\A$ is dense in $X$. Since $X$ is f.g. as a
$B$-module, we can find a finite number of generators in
$X\r\hat\A$, cf~\cite{J}. Since $\xi\r\hat\A$ is a finite
dimensional space for any $\xi\in X\r\hat\A$, we conclude that
there exists a finite dimensional non-degenerate
$\hat\A$-submodule $X_0$ of $X$ such that $X_0B=X$. There exists a
finite dimensional unitary corepresentation $V\in A\otimes B(H_V)$
such that the right $\hat\A$-modules $H_V$ and $X_0$ are
isomorphic. Fix an isomorphism $T_0\colon H_V\to X_0$ and define
$T\colon H_V\otimes B\to X$ by $T(\xi\otimes b)=(T_0\xi)b$. The
mapping $T$ is a surjective morphism of $B$-modules. Since
$H_V\otimes B$ is a f.g. Hilbert $B$-module, it makes sense to
consider the polar decomposition of $T^*$, $T^*=u|T^*|$. Then
$|T^*|$ is an invertible morphism of the $B$-module $X$, and
$u\colon X\to H_V\otimes B$ is a $B$-module mapping such that
$u^*u=\iota$. Property (ii) in the definition of
$\alpha$-equivariant Hilbert modules together with non-degeneracy
ensure that $T^*$ is a morphism of $\alpha$-equivariant modules,
hence so are $|T^*|$, $u=T^*|T^*|^{-1}$ and $u^*$. Thus $u$ is an
isomorphism of the $\alpha$-equivariant $B$-module $X$ onto
$H_{V,p}$ with $p=uu^*$. \ep

Thus $K^{\hat\A}_0(B)$ can equivalently be described as the
Grothendieck group of the semigroup of isomorphism classes of
$\alpha$-equivariant f.g. Hilbert $B$-modules. Note that $K^{\hat
\A}_1(B)$ as defined in \cite{BS} is identical to $K_1$ of the
Banach category of $\alpha$-equivariant f.g. right Hilbert
$B$-modules \cite{K}. To get the pairing with cyclic cohomology,
one needs a continuous version of equivariant cyclic cohomology,
but we shall not discuss this any further here.

Consider the C$^*$-algebra crossed product $B\rtimes\hat A$, that
is, the C$^*$-algebra generated by $\alpha(B)(\hat A\otimes 1)$,
where we consider $A$ and $\hat A$ as subalgebras of
$B(L^2(A,\varphi))$. Note that as $\hat A$ is non-unital, the
algebra $B$ should be considered as a subalgebra of
$M(B\rtimes\hat A)$. As before, we write just $b\omega$ instead of
$\alpha(b)(\omega\otimes1)$. This notation agrees with the
notation introduced earlier in the sense that $\B\rtimes\hat\A$
becomes a dense $*$-subalgebra of $B\rtimes\hat A$. In the proof
below we denote by $\B_s$ the spectral subspace of $B$
corresponding to an equivalence class $s$ of irreducible unitary
corepresentations of $(A,\D)$, so
$$
\B_s=B\r B(H_s)=B\r I_s,
$$
where we identify $\hat\A$ with $\oplus_{t\in I}B(H_t)$ and where
$I_s$ is the unit in $B(H_s)$.

\begin{lem} \label{4.3}
We have $(\B\rtimes\hat\A)(B\rtimes\hat A)(\B\rtimes\hat\A)\subset
\B\rtimes\hat\A$.
\end{lem}

\bp Let $p$ be a central projection in $\hat\A$. It is enough to
prove that $p(B\rtimes\hat A)p\subset \B\rtimes\hat\A$. Suppose a
net $\{x_n\}_n$ in $\B\rtimes\hat\A$ with $px_np=x_n$
converges to an element $x\in B\rtimes\hat A$. Fix a basis
$\omega_1,\dots,\omega_m$ in $\hat\A p$. Then
$x_n=\sum_ib_i(n)\omega_i$ for some $b_i(n)\in\B$. Since
$$
x_n=px_n=\sum_i(b_i(n)\r \hat S^{-1}(p_{(1)}))p_{(0)}\omega_i,
$$
we may assume that there exists a finite set $F$ of equivalence
classes of irreducible unitary corepresentations of $(A,\D)$ such
that $b_i(n)\in\oplus_{s\in F}\B_s$ for all $n$ and $i$. Namely,
if $p=\sum_{t\in F_0}I_t$, then for $F$ we can take the set of
irreducible components of $\bar t_1\times t_2$ for $t_1,t_2\in
F_0$. Fix matrix units $m^s_{jk}$ in $B(H_s)$, and let
$u^s_{jk}\in\A$ be the corresponding matrix coefficients. Then
there exists uniquely defined elements
$b^s_i(n),b^s_{ijk}(n)\in\B_s$ such that $b_i(n)=\sum_{s\in
F}b^s_i(n)$ and
$$
\alpha(b^s_i(n))=\sum_{j,k}u^s_{jk}\otimes b^s_{ijk}(n).
$$
By assumption the net
$$
\sum_i\alpha(b_i(n))(\omega_i\otimes1)=\sum_{s\in
F}\sum^m_{i=1}\sum_{j,k}u^s_{jk}\omega_i\otimes b^s_{ijk}(n)
$$
converges in $K(L^2(A,\varphi))\otimes B$. Since the map
$\A\otimes\hat\A\to K(L^2(A,\varphi))$, $a\otimes\omega\mapsto
a\omega$, is injective, we conclude that the nets
$\{b^s_{ijk}(n)\}_n$ converge in $B$ for all $s,i,j,k$. Then
$b^s_i(n)=\sum_{j,k}\eps(u^s_{jk})b^s_{ijk}(n)$ converge to
elements $b^s_i\in\B_s$, so $x=\sum_s\sum_i b^s_i\omega_i$ is
an element of $\B\rtimes\hat\A$. \ep

Let $V$ be a finite dimensional unitary corepresentation of
$(A,\D)$ and $p\in B(H_V)\otimes\B$ a $\hat\A$-invariant
projection. As in Subsection~\ref{s3.1} we can consider
$X_p=p(H_V\otimes\B)$ as a $(\B\rtimes\hat\A)$-module. The same
argument as we used there shows that it is isomorphic to a f.g.
projective $(\B\rtimes\hat\A)$-module $q(\B\rtimes\hat\A)^n$ for
some idempotent $q\in\Mat_n(\B\rtimes\hat\A)$. More precisely,
if $V=V^s$ is
irreducible and $p=1$, then for $q$ we can take any
one-dimensional idempotent in $B(H_{\bar s})\subset\hat\A$.
Conversely, if $q\in\Mat_n(\B\rtimes\hat\A)$, then there exists an
idempotent $e$ in $\hat\A$ such that $q(1\otimes e)=q$ in
$\Mat_n(\7C)\otimes(\B\rtimes\hat\A)$.
Then the $(\B\rtimes\hat\A)$-module $q(\B\rtimes\hat\A)^n$ is
a direct summand of the module $(e(\B\rtimes\hat\A))^n
\cong (e\hat\A\otimes\B)^n$, so it is of the form $X_p$ for some
$V$ and $p$. We can now consider $q$ as an idempotent in
$\Mat_n(B\rtimes\hat A)$, so the corresponding projective
$(B\rtimes\hat A)$-module is just
$X_p\otimes_{\B\rtimes\hat\A}(B\rtimes\hat A)$. Given two
$\hat\A$-invariant projections $p\in B(H_V)\otimes\B$ and $p'\in
B(H_{V'})\otimes\B$, it is clear that they are
equivalent if and only if the $\alpha$-equivariant modules
$H_{V,p}$ and $H_{V',p'}$ are isomorphic, if and only if the
$(\B\rtimes\hat\A)$-modules $X_p$ and $X_{p'}$ are isomorphic, if
and only if the corresponding idempotents $q$ and $q'$ in
$\cup_{n\in\7N}\Mat_n(\B\rtimes\hat\A)$ are equivalent. But
according to Lemma~\ref{4.3} equivalence of idempotents in
$\cup_{n\in\7N}\Mat_n(\B\rtimes\hat\A)$ is the same as equivalence
of idempotents in $\cup_{n\in\7N}\Mat_n(B\rtimes\hat A)$. We
summarize this discussion in the following theorem.

\begin{thm} \label{Julg2}
There are bijective correspondences between \enu{i} the
isomorphism classes of $\alpha$-equivariant f.g. right Hilbert
$B$-modules; \enu{ii} the equivalence classes of
$\hat\A$-invariant projections in $B(H_V)\otimes B$ for finite
dimensional unitary corepresentations $V$ of $(A,\D)$; \enu{iii}
the equivalence classes of idempotents in
$\cup_{n\in\7N}\Mat_n(\B\rtimes\hat\A)$; \enu{iv} the equivalence
classes of projections in $K(l^2(\7N))\otimes(B\rtimes\hat A)$.

In particular, $K^{\hat \A}_0(B)\cong K^{\hat \A}_0(\B)\cong
K_0(\B\rtimes\hat\A)\cong K_0(B\rtimes\hat A)$.
\end{thm}

\bp  The bijective correspondences are already explained. It
remains to note that though the algebras $\B\rtimes\hat\A$ and
$B\rtimes\hat A$ are in general non-unital, we can choose an
approximate unit $\{e_n\}_n$ in $B\rtimes\hat A$ consisting of
central projections in $\hat\A$, so that $K_0(B\rtimes\hat A)$ and
$K_0(\B\rtimes\hat\A)$ can be described in terms of idempotents in
$\cup_{m\in\7N}\cup_n\Mat_m(e_n(B\rtimes\hat A)e_n)
=\cup_{m\in\7N}\cup_n\Mat_m(e_n(\B\rtimes\hat\A)e_n)
=\cup_{m\in\7N}\Mat_m(\B\rtimes\hat\A)$. \epp

\subsection{Equivariant vector bundles on homogeneous spaces}

In this section we consider quantum homogeneous spaces and prove
an analogue of the geometric fact that an equivariant vector
bundle on a homogeneous space is completely determined by the
representation of the stabilizer of some point in the fiber over
this point. It seems that in the quantum case such a result is
less obvious than in the classical case, and we deduce it from
Theorem~\ref{Julg2} combined with the Takesaki-Takai duality.

We consider a closed subgroup of $(A,\D)$, that is, a compact
quantum group $(A_0,\D_0)$ together with a surjective
$*$-homomorphism $P\colon A\to A_0$ such that $(P\otimes
P)\D=\D_0P$. Consider the right coaction $\D_R\colon A\to A\otimes
A_0$ of $(A_0,\D_0)$ on $A$ given by $\D_R(a)=(\iota\otimes
P)\D(a)$, and the corresponding quotient space $B=A^{\D_R}=\{a\in
A\,|\,\D_R(a)=a\otimes1\}$ together with the left coaction
$\alpha=\D|_B\colon B\to A\otimes B$ of $(A,\D)$ on $B$. Denote by
$\Ir0$ the set of equivalence classes of irreducible unitary
corepresentations of $(A_0,\D_0)$. For each $t\in\Ir0$ fix a
representative $V^t\in A_0\otimes B(H_t)$ of this class. On
$H_t\otimes A$ define a right $A_0$-comodule structure by
$$
\delta_t\colon H_t\otimes A\to H_t\otimes A\otimes A_0,\ \
\delta_t(\xi\otimes
a)=V^t_{31}(\xi\otimes1\otimes1)(1\otimes\D_R(a)).
$$
The comultiplication $\D$ induces also a left $A$-comodule
structure on $H_t\otimes A$ by
$$
\delta\colon H_t\otimes A\to A\otimes H_t\otimes A, \ \
\delta(\xi\otimes a)=(1\otimes\xi\otimes1)\D(a)_{13}.
$$
Set
$$
X_t=(H_t\otimes A)^{\delta_t}=\{x\in H_t\otimes
A\,|\,\delta_t(x)=x\otimes1\}.
$$
Since $(\delta\otimes\iota)\delta_t=(\iota\otimes\delta_t)\delta$,
the projection $E_t\colon H_t\otimes A\to X_t$,
$E_t=(\iota\otimes\varphi_0)\delta_t$, has the property
$\delta E_t=(\iota\otimes E_t)\delta$, so $\delta$
induces a structure of an $\alpha$-equivariant right $B$-module on
$X_t$. Here $\varphi_0$ denotes the Haar state on $A_0$. To
introduce a $B$-valued inner product on it, consider the
conditional expectation $E\colon A\to B$ given by
$E=(\iota\otimes\varphi_0)\D_R$. Then set
$$
\<\xi_1\otimes a_1,\xi_2\otimes a_2\>=(\xi_2,\xi_1)E(a^*_1a_2).
$$

\begin{lem} \label{4.5}
The module $X_t$ is an $\alpha$-equivariant f.g. right Hilbert
$B$-module.
\end{lem}

\bp Consider the left action of $\hat\A_0$ on $A$ given by
$\omega\l a=(\iota\otimes\omega)\D_R(a)$. Let $A_t=B(H_{\bar t})\l
A$ be the spectral subspace of $A$ corresponding to $t$.
Equivalently $A_t$ can be described as the linear span of elements
$(\iota\otimes\varphi_0(u\cdot))\D_R(a)$, where $a\in A$ and $u$
is a matrix coefficient of $V^t$. Since $X_t=E_t(H_t\otimes A)$,
it follows that $X_t\subset H_t\otimes A_t$. Thus it is enough to
prove that $A_t$ considered as a $B$-module with inner product
$\<a_1,a_2\>=E(a^*_1a_2)$ is a  f.g. right Hilbert $B$-module.

To check that $A_t$ is complete, choose matrix units $e_{ij}$ in
$B(H_{\bar t})$ and consider the corresponding matrix coefficients
$v_{ij}$ for $V^{\bar t}$. Then define  $B$-module mappings
$E_{ij}\colon A_t\to A_t$ by letting $E_{ij}(a)=e_{ij}\l a$, so
$$
\D_R(a)=\sum_{i,j}E_{ij}(a)\otimes v_{ij}.
$$
The orthogonality relations imply that the $B$-valued inner
product $\<a_1,a_2\>'=\sum_{i,j}E_{ij}(a^*_1)E_{ij}(a_2)$ defines
an equivalent norm. Since $a=\sum_{i,j}\eps_0(v_{ij})E_{ij}(a)$
for $a\in A_t$, where $\eps_0$ is the counit on $\A_0$, it follows
that the topology on $A_t$ is the usual norm topology inherited
from $A$. Thus $A_t$ is complete.

It remains to show that $A_t$ is finitely generated. This is, in
fact, known (see e.g.~\cite{M} for a more general statement), but
we include a proof for the sake of completeness. Since
$A^*_t=A_{\bar t}$, it is enough to show that $A_t$ is finitely
generated as a left $B$-module. Let $U\in A\otimes B(H)$ be a
finite dimensional unitary corepresentation of $(A,\D)$. Then
$(P\otimes\iota)(U)$ is a unitary corepresentation of
$(A_0,\D_0)$. We can decompose it into irreducible ones, and then
the image of the space of matrix coefficients of $U$ under $P$ is
precisely the space of matrix coefficients of those irreducible
unitary corepresentations. Since $P(\A)$ is dense in $A_0$, using
orthogonality relations, we conclude that there exists an
irreducible unitary corepresentation $U$ such that $V^{\bar t}$ is
a subcorepresentation of $(P\otimes\iota)(U)$. Thus we can
identify $H_{\bar t}$ with a subspace of $H$ and complete the
orthonormal basis $\xi_1,\dots,\xi_m$ in $H_{\bar t}$ to an
orthonormal basis $\xi_1,\dots,\xi_n$ in $H$. Consider the
corresponding matrix coefficients $u_{ij}$, $1\le i,j\le n$, of
$U$. Then $P(u_{ij})=v_{ij}$ for $1\le i,j\le m$, so the elements
$b_{ij}=u_{ij}$, $1\le i\le n$, $1\le j\le m$, lie in $A_t$. We
claim that they generate $A_t$ as a left $B$-module. Indeed, for
$a\in A_t$ set
$$
a_{ij}=\sum^m_{k=1}E_{kj}(a)b^*_{ik}, \ \ 1\le i\le n, \
1\le j\le m.
$$
Since $\D_R(E_{kj}(a))
=\sum^m_{p=1}E_{pj}(a)\otimes v_{pk}$, we get
$$
\D_R(a_{ij})=\sum^m_{k,p,q=1}(E_{pj}(a)\otimes v_{pk})
(u^*_{iq}\otimes v^*_{qk})
=\sum^m_{p=1}E_{pj}(a)u^*_{ip}\otimes1=a_{ij}\otimes1,
$$
so $a_{ij}\in B$. Then we compute
$$
\sum^n_{i=1}\sum^m_{j=1}a_{ij}b_{ij}=\sum^n_{i=1}\sum^m_{j,k=1}
E_{kj}(a)u^*_{ik}u_{ij}=\sum^m_{j=1}E_{jj}(a)=I_{\bar t}\l a=a.
$$
\epp

Now we can formulate the main result of this subsection.

\begin{thm} \label{bundles}
Any $\alpha$-equivariant f.g. right Hilbert $B$-module is
isomorphic to a unique finite direct sum of modules $X_t$,
$t\in\Ir0$.

In particular, $K^{\hat\A}_0(B)$ is a free abelian group with basis
$[X_t]$, $t\in\Ir0$.
\end{thm}

To prove Theorem we first realize $B\rtimes\hat A$ as a subalgebra
of $B(\LL)$. Since $B\subset A$, we can consider $B\rtimes\hat A$
as a subalgebra of $A\rtimes\hat A$. Then we can apply the
following particular case of the Takesaki-Takai duality, see
e.g.~\cite{BS2}.

\begin{lem}
The representations of $A$ and of $\hat A$ on $\LL$ define an
isomorphism of $A\rtimes\hat A$ onto $K(\LL)$.

Specifically, the homomorphism $B(\LL)\to B(\LL\otimes\LL)$,
$x\mapsto V(x\otimes 1)V^*$, sends $\omega\in\hat A$ to
$\omega\otimes1$ and $a\in A$ to $\D(a)$, where $V=(J\hat
J\otimes1)\Sigma W\Sigma(J\hat J\otimes1)$, $J$ is the modular
involution corresponding to $\varphi$, $\hat J$ is the modular
involution corresponding to the right invariant Haar weight on
$\hat A$, $\Sigma$ is the flip, and $W$ is the multiplicative
unitary for $(A,\D)$.\epp
\end{lem}

Thus there is no ambiguity in writing $a\omega$ irrespectively of
whether we consider it as an element of $A\rtimes\hat A$ or as an
operator on $\LL$.

Define a unitary $U_0$ on $\LL\otimes L^2(A_0,\varphi_0)$ by
$$
U_0(a\xi_\varphi\otimes\xi)=\D_R(a)(\xi_\varphi\otimes\xi).
$$
Then define a right coaction $\tilde\D_R\colon K(\LL)\to
M(K(\LL)\otimes A_0)$ of $(A_0,\D_0)$ on $K(\LL)$ by setting
$\tilde\D_R(x)=U_0(x\otimes1)U^*_0$. We obviously have
$\tilde\D_R(a)=\D_R(a)$ for $a\in A$. Since the representation of
$\hat\A$ on $\LL$ by definition has the property $\omega
a\xi_\varphi=(\hat S^{-1}(\omega)\otimes1)\D(a)\xi_\varphi$, we
see also that $\tilde\D_R(\omega)=\omega\otimes1$. It follows that
we can identify $B\rtimes\hat A$ with $(A\rtimes\hat
A)^{\tilde\D_R}$.

Choose an irreducible unitary corepresentation $U\in A\otimes
B(H)$ such that $V^t$ is a subcorepresentation of
$(P\otimes\iota)(U)$, and choose an orthonormal basis
$\xi_1,\dots,\xi_n$ in $H$ such that $\xi_1,\dots,\xi_m$ is an
orthonormal basis in $H_t$. Let $u_{ij}$, $1\le i,j\le n$, be the
corresponding matrix coefficients for $U$, and $v_{ij}$, $1\le
i,j\le m$, the matrix coefficients for $V^t$. Fix $i$, $1\le\ i\le
n$.

\begin{lem} \label{4.8}
The map $H_t\otimes\A\to K(\LL)$, $\xi_j\otimes a\mapsto
u_{ij}I_0a$, defines a right $(\A\rtimes\hat\A)$-module and right
$\A_0$-comodule isomorphism of $H_t\otimes\A$ onto
$p_t(\A\rtimes\hat\A)$, where $p_t\in K(\LL)$ is the projection
onto the space spanned by $u_{ij}\xi_\varphi$, $1\le j\le m$.
\end{lem}

\bp We have already used the isomorphisms $\A\cong
I_0\hat\A\otimes\A\cong I_0(\A\rtimes\hat\A)$, where
$I_0\in\hat\A$ is the projection such that
$I_0\omega=\hat\eps(\omega)I_0$. In other words, the map $\A\to
I_0(\A\rtimes\hat\A)$, $a\mapsto I_0a$, is an
$(\A\rtimes\hat\A)$-module isomorphism. Note also that since $I_0$
is the projection onto $\7C\xi_\varphi$, the operator $u_{ij}I_0$
is, up to a scalar, a partial isometry with initial space
$\7C\xi_\varphi$ and range space $\7Cu_{ij}\xi_\varphi$. It
follows that the map in the formulation of Lemma is an injective
morphism of $(\A\rtimes\hat\A)$-modules. Moreover, since $p_t$ is
a linear combination of $u_{ij}I_0u^*_{ij}$, we have
$u_{ij}I_0(\A\rtimes\hat\A)\subset p_t(\A\rtimes\hat\A)$ and
$p_t(\A\rtimes\hat\A)\subset\sum_ju_{ij}I_0(\A\rtimes\hat\A)$.
Thus the map is an isomorphism of the $(\A\rtimes\hat\A)$-modules
$H_t\otimes\A$ and $p_t(\A\rtimes\hat\A)$. It remains to check
that it is an $\A_0$-comodule morphism. By definition we have
$$
\delta_t(\xi_j\otimes a)=V^t_{31}(\xi_j\otimes1\otimes1)
(1\otimes\D_R(a))
=\sum_k(\xi_k\otimes1\otimes v_{kj})(1\otimes\D_R(a)).
$$
On the other hand,
$$
\tilde\D_R(u_{ij}I_0a)=\sum_k(u_{ik}\otimes v_{kj})
(I_0\otimes1)\D_R(a).
$$
Thus the map is indeed an $\A_0$-comodule morphism.
\ep

\bpp{Theorem~\ref{bundles}} Consider the
$(\B\rtimes\hat\A)$-module $\chi_t=(H_t\otimes\A)^{\delta_t}$. As
an immediate consequence of Lemma~\ref{4.8} we have $\chi_t\cong
p_t(\B\rtimes\hat\A)$. By Theorem~\ref{Julg2} we get a f.g.
$\alpha$-equivariant right Hilbert $B$-module $\tilde
X_t=\chi_t\otimes_\B B$. The inclusion $\chi_t\hookrightarrow X_t$
induces a morphism $\tilde X_t\to X_t$ of $\alpha$-equivariant
f.g. right Hilbert $B$-modules. This map is surjective, since by
the proof of Lemma~\ref{4.5}, we know that $\chi_t B=X_t$. The map
is also injective, since it is equivariant and injective on
$\chi_t$, and since any equivariant $B$-module map on $\tilde X_t$
maps $\chi_t$ into itself. Thus $\tilde X_t\cong X_t$. By
Theorem~\ref{Julg2} it is thus enough to show that any projection
in $K(l^2(\7N))\otimes(B\rtimes\hat A)$ is equivalent to a unique
direct sum of projections $p_t$, $t\in\Ir0$.

Note that $(U_0)_{21}\in M(A_0\otimes K(\LL))$ is a unitary
corepresentation of $(A_0,\Delta_0)$. Thus we can decompose the
space $\LL=\oplus_t(H_t\otimes L_t)$ so that
$(U_0)_{21}=\oplus_t(V_t\otimes 1)$. Then $B\rtimes\hat
A\cong\oplus_t K(L_t)$, so to prove Theorem, it is enough to show
that under this isomorphism $p_t$ becomes a minimal projection in
$K(L_t)$. Equivalently, we must show that the restriction of
$(U_0)_{21}$ to $p_t\LL$ is isomorphic to $V^t$. This is indeed
the case since in the notation of the proof of Lemma~\ref{4.8} we
have
$$
(U_0)_{21}(\xi\otimes u_{ij}\xi_\varphi)
=\sum_kv_{kj}\xi\otimes u_{ik}\xi_\varphi.
$$
\epp

\bigskip

\flushleft
{Sergey Neshveyev\\
Mathematics Institute\\
University of Oslo\\
PB 1053 Blindern\\
Oslo 0316\\
Norway\\
{\it e-mail}: neshveyev@hotmail.com}

\flushleft
{Lars Tuset\\
Faculty of Engineering\\
Oslo University College\\
Cort Adelers st. 30\\
Oslo 0254\\
Norway\\
{\it e-mail}: Lars.Tuset@iu.hio.no}


\bigskip\bigskip

\begin{thebibliography}{99}

\bibitem[AK1]{AK1}
Akbarpour R., Khalkhali M., {\it Equivariant cyclic cohomology of
Hopf module algebras}, preprint math.KT/0009236.


\bibitem[AK2]{AK2}
Akbarpour R., Khalkhali M., {\it Hopf algebra equivariant cyclic
homology and cyclic homology of crossed product algebras},
preprint math.KT/0011248.

\bibitem[BS1]{BS}
Baaj S., Skandalis G., {\it $C\sp *$-alg\`ebres de Hopf et
th{\'e}orie de Kasparov {\'e}quivariante}, $K$-Theory {\bf 2}
(1989), 683--721.

\bibitem[BS2]{BS2}
Baaj S., Skandalis G., {\it Unitaires multiplicatifs et
dualit{\'e} pour les produits crois{\'e}s de $C\sp *$-alg\`ebres},
Ann. Sci. École Norm. Sup. (4){\bf 26} (1993), 425--488.

\bibitem[C]{Co}
Connes A. Noncommutative geometry.
Academic Press, Inc., San Diego, CA, 1994.


\bibitem[D]{M}
Durdevi\'c M., {\it Geometry of quantum principal bundles. II},
Rev. Math. Phys. {\bf 9} (1997), 531--607.

\bibitem[J]{J}
Julg P., {\it $K$-th{\'e}orie {\'e}quivariante et produits
crois{\'e}s}, C. R. Acad. Sci. Paris S{\'e}r. I Math. {\bf 292}
(1981), 629--632.

\bibitem[K]{K}
Karoubi M. $K$-theory. An introduction. Grundlehren der
Mathematischen Wissenschaften, Band 226. Springer-Verlag,
Berlin-New York, 1978.

\bibitem[KKL]{KKL}
Klimek S., Kondracki W., Lesniewski A., {\it Equivariant entire
cyclic cohomology. I. Finite groups}, $K$-Theory {\bf 4} (1991),
201--218.

\bibitem[L]{L}
Loday J.-L. Cyclic homology. Springer-Verlag, Berlin, 1998.

\bibitem[KMT]{KMT}
Kustermans J., Murphy G. J., Tuset L., {\it Differential calculi
over quantum groups and twisted cyclic cocycles}, J. Geom. Phys.
{\bf 44} (2003), 570--594.

\bibitem[NT]{NT}
Neshveyev S., Tuset L., {\it The Martin boundary of a discrete
quantum group}, preprint math.OA/0209270.

\end{thebibliography}
\end{document}